\documentclass [11pt] {article}
\usepackage{amsfonts}
\usepackage{amssymb}
\newcommand{\qed} {\hspace {0.1in} \rule {1.5mm} {3.5mm}}

\newtheorem{lemma}{Lemma}[section]
\newtheorem{corollary}{Corollary}[section]
\newtheorem{statement}{Statement}
\newtheorem{theorem}{Theorem}
\newtheorem{example}{Example}
\newtheorem{question}{Question}

\newtheorem{conjecture}{Conjecture}
\newtheorem{proposition}{Proposition}[section]
\newtheorem{definition}{Definition}[section]
\def\limn{\lim_{n\to\infty}}

\def\e{\epsilon}

\def\deg{\mbox{deg}\,}

\def\tr{\mbox{Tr}\,}

\def\trg{\tr_{\bG}}

\def\dim{{\rm dim}}

\def\<{\langle}
\def\>{\rangle}

\def\proof{\smallskip\noindent{\it Proof.} }

\def\bG{{\bf G}}
\def\bA{{\bf A}}
\def\bB{{\bf B}}
\def\bH{{\bf H}}

\def\bR{{\mathbb R}}
\def\bN{{\mathbb N}}
\def\bC{{\mathbb C}}

\def\deg{\mbox{deg}\,}

\def\cA{\mbox{$\cal A$}}

\def\cA{\mbox{$\cal A$}}

\def\cP{\mbox{$\cal P$}}

\def\cN{\mbox{$\cal N$}}

\def\ds{\delta_s}
\def\dr{\delta_\rho}
\def\Fo{F$\mbox{\o}$lner}

\def\rk{\mbox{Rank}}
\def\ah{\hat{A}}
\def\ha{\ah}
\begin{document}

\title{$L^2$-spectral invariants and convergent sequences of finite graphs}
\author{{\sc G\'abor Elek}
\footnote {The Alfred Renyi Mathematical Institute of
the Hungarian Academy of Sciences, P.O. Box 127, H-1364 Budapest, Hungary.
email:elek@renyi.hu, Supported by OTKA Grants T 049841 and T 037846}}
\date{}
\maketitle \vskip 0.2in \noindent
\vskip 0.2in \noindent{\bf Abstract.} Using the spectral theory 
of weakly convergent sequences of finite graphs, we prove the
uniform existence of the integrated density of states for a large class 
of infinite graphs. 

\vskip 0.3in

\newpage
\tableofcontents
\newpage
\section{Introduction}
The goal of the 2006 Oberwolfach-Mini-Workshop ``$L^2$-spectral invariants and
the integrated density of states''  was to unify the point of views
and approaches in certain areas of geometry and mathematical
physics. The aim of our paper is to make the connection between
those fields even more explicit. Let us start with a very brief
introduction to the theory of integrated density of states. 
\subsection{Laplace operators on infinite graphs and their
integrated densities of states}
Let $G$ be an infinite connected graph with bounded vertex degrees. We
say that $G$ is {\bf amenable} if there exists a sequence of finite
connected spanned subgraphs $\{G_n\}^\infty_{n=1}$ such that
$$\limn \frac{|\partial G_n|}{|V(G_n)|}=0\,,$$
where 
$$\partial G_n=\{x\in V(G_n)\,\mid \,\mbox{there exists $y\in V(G)\backslash
V(G_n)$ such that $(x,y)\in E(G_n)$}\}\,.$$
Such sequences of subgraphs are called {\bf \Fo-sequences}. Now let
$\Delta_n:L^2(V(G_n))\to L^2(V(G_n))$ be the Laplacian operator
$$\Delta_nf(x)=\deg(x) f(x)-\sum_{\{y\,\mid (x,y)\in E(G_n)\}} f(y)\,,$$
where the degree of $x$ is considered in the subgraph $G_n$. Then $\Delta_n$
is a finite dimensional positive self-adjoint operator.
Let 
$$N_{\Delta_n}(\lambda):=\frac{|\{\mbox{eigenvalues of $\Delta_n$ not larger
than $\lambda$ (with multiplicities)}\}|}{|V(G_n)|}\,.$$
We call $N_{\Delta_n}$ the {\bf normalized spectral distribution function}
of $\Delta_n$. We say that the {\bf integrated density of states} 
exists for the Laplacian of the graph $G$ if there exists a right continuous
monotone function $\sigma$ such that for any \Fo-sequence
$\{G_n\}^\infty_{n=1}$:
$$\limn N_{\Delta_n}(\lambda)=\sigma(\lambda)\,,$$
if $\lambda$ is a continuity point of $\sigma$.
We say that the integrated density of state uniformly exists if
$\{N_{\Delta_n}\}^\infty_{n=1}$ uniformly converge to $\sigma$.
\begin{question}
For which amenable graphs $G$ does the integrated density of states exist ?
\end{question}
Let $H$ be the $2$-dimensional lattice and $K$ be the $3$-dimensional
lattice. Construct a new graph $G$ by identifying a vertex of $H$ with
a vertex of $K$.
Then the integrated density of states clearly does not exist for the Laplacian
of $G$. This example
suggests that one needs some sort of homogeneity in the local geometry of $G$.
\subsection{The periodic case}
Let $\Gamma$ be a countable group and $L^2(\Gamma)$ be the
Hilbert-space of the formal sums $\sum_{\gamma\in\Gamma}
a_\gamma\cdot\gamma$, where $a_\gamma\in\bC$ and
$\sum_{\gamma\in\Gamma} |a_\gamma|^2<\infty$. Notice that $\Gamma$
unitarily acts on $L^2(\Gamma)$ by

\noindent
 $L_\delta(\sum_{\gamma\in\Gamma} a_\gamma\cdot\gamma)=
\sum_{\gamma\in\Gamma} a_{\delta^{-1}\gamma}\cdot\gamma\,.$ Hence
one can represent the complex group algebra $\bC\Gamma$ as bounded
operators by left convolutions. The weak closure of  $\bC\Gamma$ in
$B(L^2(\Gamma))$ is the group von Neumann algebra $\cN\Gamma$. The
group von Neumann algebra has a natural trace:
$$\tr_\Gamma(A)=\langle A(1),1\rangle\,,$$
where $1\in L^2(\Gamma)$ is identified with the unit element of the
group. Let $B\in\cN\Gamma$ be a self-adjoint element, then
by the spectral theorem of von Neumann
$$B=\int_{-\infty}^{\infty}\lambda\,d E^B_\lambda\,,$$
where $E^B_\lambda=\chi_{[-\infty,\lambda]}(B)\,.$ We can associate
a {\it spectral measure} $\mu_B$ to our operator $B$ by
$$\mu_B[-\infty,\lambda]=\sigma_B(\lambda)=\tr_\Gamma
E^B_\lambda\,.$$ Note that the jumps of $\sigma_B$ are associated to
the eigenspaces of $B$. 
\noindent
Now let $\Gamma$ be a finitely generated amenable group and
$\mbox{Cay}(\Gamma,S)$ be the Cayley-graph of $\Gamma$ with respect
to a symmetric generating set $S$. Then the Laplacian of $G$ can be
regarded as the element
$$\Delta_G=\|S\|{\bf 1}-\sum_{s\in S} s\in \bC\Gamma\,.$$
Hence the obvious candidate for the integrated density of states is
the spectral measure $\sigma_{\Delta_G}$. In fact one has the following
result.
\begin{statement} \label{s1} \cite{DLM} For the Cayley-graph of an amenable
group and a \Fo-subgraph sequence $\{G_n\}^\infty_{n=1}$ if
$B\in \bC(\Gamma)$ is a self-adjoint element then 
$\{N_{B_n}\|^\infty_{n=1}$ uniformly converges to $\sigma_B$, where
$B_n=p_n B i_n$ and $p_n:L^2(\Gamma)\to L^2(G_n)$
is the natural projection operator, $i_{F_n}:L^2(F_n)\to
L^2(\Gamma)$ is the adjoint of $p_{F_n}$ the natural imbedding
operator.
\end{statement}

Similar approximation theorem holds for residually finite groups
in the weaker sense.
Let $\Gamma$ be a finitely generated residually finite group with
finite index normal subgroups
$$\Gamma\rhd N_1 \rhd N_2\rhd\dots, \cap^\infty_{k=1} N_k=\{1\}\,.$$
\begin{statement} \label{Luck} \cite{Lueck1}
Let $B\in\bC\Gamma$ be a self-adjoint element and let
$\pi_k(B)=B_k\in\bC(\Gamma/N_k)$ be the associated finite
dimensional linear operators, where $\pi_k:\Gamma\to \Gamma/N_k$ are
the quotient maps. Then  the spectral distribution functions $N_{B_k}$
converge at any continuity point of $\sigma_B$.
\end{statement}
According to the Strong Approximation Conjecture of L\"uck the convergence
in Statement \ref{Luck} is always uniform. Note that the conjecture holds
for
amenable groups \cite{Elekstrong}.

\subsection{The aperiodic case}
In \cite{LS} (see also \cite{Lenz} for a short exposition) 
Lenz and Stollmann studied 
graphs constructed by Delone sets in $\bR^n$. 
In these graphs each neighborhood pattern can be seen in a given
frequency but they do not have
any sort of global symmetries. Instead of the Laplacians they considered
finite range pattern-invariant operators. These operators can be viewed
as the aperiodic analogs of the elements of the group algebra.
They proved the following result.
\begin{statement} \label{s3}
Let $G$ be a graph of a Delone-set and $\{G_n\}^\infty_{n=1}$ be
a \Fo-sequence. Also, let $A$ be self-adjoint finite range pattern-invariant
operator on $G$. Then the normalized spectral distributions 
$\{N_{A_n}\}^\infty_{n=1}$ converge uniformly to an integrated density of
states $\sigma_A$ that does not depend on the choice of the \Fo-sequence.
\end{statement}
Note that the weak convergence of
the spectral distributions was already
established by Kellendonk \cite{Kell} and by Hof \cite{Hof}. 
 The obstacle what Lenz and Stollmann
had to overcome was the possible discontinuity of the integrated density
of states due to the existence of finitely supported eigenfunctions. This
phenomenon did not occur in the case of lattices. Note however that
for certain amenable groups Grigorchuk and Zuk proved the existence
of finitely supported eigenfunctions \cite{GZ} for the Laplacian operators.
\subsection{Our results}
In Section \ref{graph} we study weakly convergent sequences of finite graphs
introduced by Benjamini and Schramm \cite{BS}. These graph sequences
are exactly the ones for which each neighborhood pattern can be seen
at a certain frequency. We introduce two notions: {\bf strong graph
convergence}
and {\bf antiexpanders}. Strong graph convergence immediately ensures
the uniform convergence of the normalized spectral distributions of
the Laplacians. The notion of antiexpander seems to be the
 right notion of amenability
in the world of weakly convergent graph sequences.
We prove that one can always pick strongly convergent subsequences
from weakly convergent antiexpander sequences and conjecture that
weakly convergent antiexpander sequences are actually always strongly
convergent. In Section \ref{algebra} we associate von Neumann algebras
to weakly convergent graph sequences via finite range pattern-invariant
operator sequences. Note that von Neumann algebra was considered in \cite{LS}
as well, using the action of $R^n$ on the tiling space generated by
the Delone-set. In this paper we do not have any group action only
the ``{\it  statistical symmetry}'' given by the existence of 
the pattern frequency.
We construct a limit operator in our von Neumann algebra for
finite range pattern-invariant operator sequences.
Generalizing Statement \ref{Luck} we prove that the
spectral distribution measure of the limit operator is the integrated
density of states. 
Then in Section \ref{infinite} we are dealing with amenable graphs
and as our main result we answer Question 1 for a rather large class of 
graphs generalizing Statement \ref{s1} and Statement \ref{s3}.

\vskip0.1in
\noindent
{\bf Main result} {\it (Theorem \ref{main}) If $G$ is an amenable graph
such that all the \Fo- sequences $\{G_n\}^\infty_{n=1}$ are weakly
convergent antiexpanders, then for any 
self-adjoint finite range pattern invariant operator $A$ the
integrated density of states $\sigma_A$ uniformly exists.}

\section{Graph sequences of bounded vertex degrees}\label{graph}
\subsection{The weak convergence of graph sequences}
First let us recall the notion of weak convergence of finite graphs
due to Benjamini and Schramm in a slightly more general form as in
\cite{BS}. A {\bf rooted
  $(d,r,X,S)$}-graph is a finite simple connected graph $G$
\begin{itemize}
\item with a distinguished vertex $x$ (the root),
\item such that $\deg(y)\leq d$ for any $y\in V(G)$,
\item such that $d_G(x,z)\leq r$ for any $z\in V(G)$, where
$d_G$ denotes the usual shortest path distance,
\item the vertices of $G$ are colored by the elements of the set $X$,
\item the directed edges of $G$ are colored by the elements
of the set $S$ (that is both $\overrightarrow{(x,y)}$
and $\overrightarrow{(x,y)}$, if $(x,y)\in E(G)$).
\end{itemize}
In general, we shall call finite graphs with vertices colored by $X$ and 
directed edges colored by $S$; $(X,S)$-graphs. Two rooted $(d,r,X,S)$-graphs
$G$ and $H$ are 
called {\bf rooted isomorphic} if there exists a graph isomorphism
between them mapping root to root, preserving both the vertex-colorings
 and the edge-colorings.

\noindent
Let $\cA(d,r,X,S)$ denote the finite set of rooted isomorphism classes
of rooted $(d,r,X,S)$-graphs. Now let $G$ be an arbitrary finite
$(X,S)$-graph. Then for any $r\geq 1$ we can associate to $G$ a probability
distribution on $\cA(d,r,X,S)$ by
$$p_G(\alpha)=\frac{|T(G,\alpha)|}{|V(G)|},$$
where $T(G,\alpha)$ denotes the set of vertices $z\in V(G)$ such that
the $r$-neighbourhood of $z$; $B_r(z)$ represents the class $\alpha$.
\begin{definition}
Let $\bG=\{G_n\}^\infty_{n=1}$ be a 
sequence of finite connected \\$(X,S)$-graphs
such that $V(G_n)\to\infty$ as $n$ tends to $\infty$. Then we say that
$\bG=\{G_n\}^\infty_{n=1}$ is {\bf weakly convergent} if for any $r\geq 1$
and $\alpha\in \cA(d,r,X,S)$, $\limn p_{G_n}(\alpha)$ exists.
\end{definition}
Note that if $|X|=1$, $ |S|=1$ then this is the usual definition of weak
convergence for non-colored graphs.

\subsection{Strong convergence}
Let $G$ be an $(X,S)$-graph and $y\in V(G)$. Then the {\bf star} of $y$,
$S_G(y)$ is defined as follows;
\begin{itemize}
\item $V(S_G(y))$ consists of $y$ and its neighbours.
\item $E(S_G(y))$ consists of the edges between $y$ and its neighbours.
\item The coloring of $S_G(y)$ is inherited from $G$ (that is $S_G(y)$ is
an $(X,S)$-graph as well).
\end{itemize}
\begin{definition}
Let $G,H$ be two not necessarily connected graphs on the vertex set $V$. Then 
the $\delta$-distance of them is defined as
$$\delta(G,H):=\frac{|\{y\in V\,\mid\, S_G(y) \ncong S_H(y)\}|}{|V|}\,.$$
\end{definition}
Note that $\cong$ means rooted $(X,S)$-colored isomorphism not merely graph
isomorphism.
\begin{lemma}
$\delta$ defines a metric on the $(X,S)$-graphs with vertex set $V$.
\end{lemma}
\proof
Clearly, $\delta(G,H)=0$ if and only if $G=H$. Also, $\delta(G,H)=\delta(H,G)$.
Let us check the triangle inequality. Let $G,H,J$ be three graphs on $V$.
$$\{y\in V\,\mid\, S_G(y) \ncong S_J(y)\}\subseteq
\{y\in V\,\mid\, S_G(y) \ncong S_H(y)\}\cup \{y\in V\,\mid\, S_H(y) \ncong 
S_J(y)\}\,.$$
That is $\delta(G,J)\leq \delta(G,H)+\delta(H,J).$\qed

\noindent
Let $\sigma\in S(V)$ be a permutation of the vertices. Then $H^\sigma$ denotes
the $(X,S)$-graph on $V$, where $(\sigma(x),\sigma(y))\in E(H^\sigma)$ if
and only if $(x,y)\in E(H)$. Also, $\sigma(x)$ in $H^\sigma$ is colored the
same way as $x$ is colored in $H$, respectively $(\sigma(x),\sigma(y)$ is
colored the same way as $(x,y)$ is colored in the graph $H$. Hence we can
define
$$\delta_s(G,H)=\inf_{\sigma\in S(V)} \delta(G, H^\sigma)\,.$$
\begin{lemma}
$\delta_s$ defines a metric on the isometry classes of $(X,S)$-graphs
with vertex set $V$.
\end{lemma}
\proof
Clearly, $G \cong H$ if and only if $\delta_s(G,H)=0.$ Also,
$$\ds(G,H)=\inf_{\sigma\in S(V)} \delta(G,H^\sigma)=
\inf_{\sigma\in S(V)} \delta(G^{\sigma^{-1}},H)=\ds(H,G)\,.$$
Now let $G,H,J$ be three graphs with vertex set $V$ and let
$\ds(G,H)=\delta(G,H^{\sigma_1})$ and $\ds(H,J)=\delta(H,J^{\sigma_2})\,.$
Then
$\ds(H,J)=\delta(H^{\sigma_1}, J^{\sigma_2\sigma_1})$, hence
$$\ds (G,J)\leq \delta(G, J^{\sigma_2\sigma_1})\leq
\ds(G,H)+\ds(H,J)\,.\quad\qed
$$

\noindent
Now we are ready to define the geometric distance of two arbitrary
finite connected $(X,S)$-graphs. Let $G$ and $H$ be $(X,S)$-graphs. Also,
let $q,r$ be two integers such that $q|V(G)|=r|V(H)|$. Denote by $qG$
the union of $q$ disjoint copies of $G$. Then $\ds(qG,rH)$ is well-defined
since the graphs $qG$ and $rH$ can be represented on the same vertex set.
\begin{definition}
The geometric distance of the finite connected $(X,S)$-graphs $G$ and $H$
is defined as
$$\dr(G,H)=\inf_{\{q,r\,\mid\, q|V(G)|=r|V(H)|\}} \ds(qG,rH)\,.$$
\end{definition}
\begin{proposition}
$\dr$ defines a metric on the set of isomorphism classes of finite
connected $(X,S)$-graphs.
\end{proposition}
\proof
Clearly, $\dr(G,H)=\dr(H,G)$. Now let us check the triangle
inequality. Let $\dr(G,H)\geq \ds(qG,rH)+\e$ and
$\dr(H,J)\geq \ds(sH,tJ)+\e$. Obviously,
$\ds(qG,rH)\geq \ds(sqG,srH)$.
Thus
$$\dr(G,H)+\dr(H,J)\geq
\ds(sqG,srH)+\ds(srH,rtJ)+2\e\geq\ds(sqG,rtJ)+2\e\geq
\dr(G,J)+2\e\,.$$
Letting $\e\to 0$ we obtain the triangle inequality. The last step
is to prove that if $\dr(G,H)=0$ then $G\cong H$. First suppose that
$|G|<|H|$. Let $tH$ and $sG$ be represented on the same vertex set $V$.
We would like to estimate $\delta(tH,(sG)^\sigma)\,.$
Since $|G|<|H|$ at least one edge in each of the $t$ copies of $H$
connects two different components in $(sG)^\sigma$ hence
$$\delta(tH,(sG)^\sigma)\geq\frac{t}{t |V(H)|}=\frac{1}{|V(H)|}\,.$$
This shows immediately that $\dr(H,G)\geq \frac{1}{|V(H)|}>0$.

\noindent
Now let us suppose that $|G|=|H|$, but $G\ncong H$. Suppose that 
$tG$ and $tH$ are represented on the
same vertex set $V$. Again, we estimate $\delta(tH,(tG)^\sigma)\,.$
For each of the $t$ copies of $H$;
\begin{itemize}
\item there exists an edge connecting two components in $(tG)^\sigma$
\item or the vertices of the particular copy are exactly the vertices of
a copy of $G$ in $(tG)^\sigma$.
\end{itemize}
Consequently in each of the $t$ copies
there exists at least one vertex such that its star in $tH$ is not
isomorphic to its star in $(tG)^\sigma$. Hence
$\delta(tH,(tG)^\sigma)\geq \frac{1}{|V(H)|}$. That is $\dr(G,H)>0$. \qed

\subsection{Strong convergence implies weak convergence}
Let $\bG=\{G_n\}^\infty_{n=1}$ be a sequence of finite connected
$(X,S)$-graphs with vertex
degree bound $d$ such that $\{G_n\}^\infty_{n=1}$ is a Cauchy-sequence
in the $\dr$-metric. Then we say that $\bG$ is a {\bf strongly convergent
graph sequence}.
\begin{proposition}
If $\bG$ is strongly convergent then $\bG$ is weakly convergent as well.
\end{proposition}
\proof
Let $r>0$ be a natural number, $\e>0$ be a real number. It is enough
to prove that there exists $\delta>0$ depending on $r,\e$ and on the
uniform degree bound $d$ such that if $\dr(G,H)<\delta$ for the $(X,S)$-
graphs $G,H$ with vertex degree bound $d$, then
$|p_G(\alpha)-p_H(\alpha)|<\e$ for any $\alpha\in\cA(d,r,X,S)$.
Let $t$ be the maximal possible number of elements in a finite connected
graph with diameter $2r$ and vertex degree bound $d$.
Let $\delta=\frac{\e}{3t}$ and suppose that
$p_G(\alpha)-p_H(\alpha)>\e$ for some $\alpha\in\cA(d,r,X,S)$ and 
also $\dr(G,H)<\delta\,.$
Then there exists $l$ and $m$, $\frac{l} {m}= \frac{|V(H)|}{|V(G)|}$,
such that $\delta(lG,(mH)^{\sigma})<2\delta\,$ (We represent $lG$ and $mH$ on
the same vertex set). Note that $p_{lG}(\alpha)=p_G(\alpha),
p_{(mH)^\sigma}(\alpha)=p_H(\alpha)\,.$ 
Let $T_1$ (resp. $T_2$) be the set of vertices
in $lG$ (resp in $(mH)^{\sigma}$ having $r$-neighbourhood isomorphic to
$\alpha$. According to our assumption
\begin{equation}
\label{e10}
|T_1|-|T_2|>\e l |V(G)|\,.
\end{equation}
The number of vertices $x$ in $lG$ such the star of $x$ in $lG$ is not
isomorphic to its star in $(mH)^{\sigma}$ is less than $2\delta l
|V(G)|$. Denote this set by $W$. Let $T'\subseteq T_1$ be the set of vertices
$z$ such that $B_r(z)$ contains an element of $W$. Observe, that if
$y\in T_1\backslash T_1'$ then $y\in T_2$. Indeed if $y\in T_1\backslash T_1'$
then the $r$-neighborhood of $y$ as an $(X,S)$-graph in $lG$ is exactly the
$r$-neighborhood of $y$ in $(mH)^\sigma$.
For $w\in W$ let $S_w$ be the set of vertices in $lG$ such that
if $x\in S_w$ then $w\in B_r(x)$. Clearly, $|S_w|\leq t$. Therefore $|T_1'|\leq
2\delta tl|V(G)|\,.$ Since $2t\delta<\e$ we are in contradiction with
(\ref{e10})\,.\qed

\vskip0.1in
{\bf Remark: }\,
The motivation for the definition of our geometric graph distance was
the graph distance $\delta_\square$
defined in \cite{BCL}. In \cite{BCL} the authors studied
the convergent sequences of dense graphs and proved (Theorem 4.1) that
a sequence of dense graphs is convergent if and only if they form a Cauchy-
sequence in the $\delta_\square$ metric. 

\subsection{Antiexpanders}
Antiexpanders (or in other words hyperfinite graph sequences) were introduced
in \cite{Elekcost}. 

\begin{definition}
Let $\bG=\{G_n\}^\infty_{n=1}$ be a sequence
of finite, connected graphs with vertex degree bound $d$. Then $G$ is an
{\bf antiexpander} if for any $\e>0$ there exists $K_{\e}>0$ such that
for any $n\geq 1$ one can remove $\e |E(G_n)|$ edges in $G_n$ such a way that
the maximal number of vertices in a component of the remaining graph $G_n'$
is at most $K_{\e}$. \end{definition}
The simplest example for an antiexpander sequence is $\{P_n\}^\infty_{n=1}$,
where $P_n$ is a path of length $n$. The reason we call such sequence
antiexpander is that if $\bH=\{H_n\}^\infty_{n=1}$ is an expander sequence
then for some $\e>0$ by removing not more than $\e|E(G_n)|$ edges from $G_n$
then at least one of the components of the remaining graph will have
size at least $\frac{1}{2} |V(G_n)|$. Thus the notion of
 antiexpanders is indeed the
opposite of the notion of expanders.
\begin{proposition}
If $\bG=\{G_n\}^\infty_{n=1}$ is a strongly convergent sequence of vertex
degree bound $d$ then $G$ is an antiexpander sequence.
\end{proposition}
\proof
First pick a constant $\e>0$. Since $\bG$ is strongly convergent there
exists a number $n_\e$ such that $\dr(G_{n_\e},G_n)<\frac{\e} {2d}$ if
$n\geq n_\e$.

\noindent{\bf Claim:} {\it If $n\geq n_\e$ then one can remove
$\e E(G_n)$ edges from $G_n$ such that
in the remaining graph $G'_n$ the maximal number of vertices in a
component is at most $|V(G_{n_\e})|$.}

\vskip0.1in
\noindent
Indeed, let us represent $pG_{n_\e}$ and $qG_n$ on the same vertex space such
a way that
$\delta(pG_{n_\e}, qG_n)<\frac{\e} {d}$\,. Let us denote by $T_n$ the set of
vertices $x$ in $qG_n$ such that the star of $x$ in $qG_n$ is not isomorphic
to the star of $x$ in $pG_{n_\e}$.
By removing edges incident to the vertices of $T_n$, all the remaining
components shall have edges from $pG_{n_\e}$. Thus the maximal number of
vertices in a component shall be at most $|V(G_{n_\e})|$.
By our assumption, 
$$\frac{|T_n|}{q|V(G_n)|}\leq\frac{\e}{d}\,.$$
Therefore we removed at most $\e q |V(G_n)|$ edges. By the pigeon hole
principle, there exists a component of $qG_n$ from which we removed at most
$\e|V(G_n)|$ edges. Thus our claim and hence the proposition itself
follows. \qed

\subsection{Weakly convergent antiexpanders}

In the previous subsections we proved that if $\bG=\{G_n\}^\infty_{n=1}$
is a strongly convergent sequence of graphs with vertex degree bound $d$,
then $\bG$ is a weakly convergent antiexpander system. The following
proposition states that at least a partial converse holds.
\begin{proposition}\label{choose}
Let $\bG=\{G_n\}^\infty_{n=1}$ be a weakly convergent 
antiexpander sequence of $(X,S)$-graphs
with vertex degree bound $d$. Then there exists a subsequence
$\{G_{n_k}\}^\infty_{k=1}$ that converges strongly as a sequence of
$(X,S)$-graphs.
\end{proposition}
\proof
Clearly, it is enough to prove that for any $\e>0$ there exists
a subsequence $\{G_{n_k}\}^\infty_{k=1}$ such that
$\dr(G_{n_l}, G_{n_m})\leq \e$ for any pair $l,m\geq 1$.
Let $\delta_1=\frac{\e}{100d}$. For any $n\geq 1$ we remove $\delta_1|E(G_n)|$
edges from $G_n$ such that in the remaining graphs $\{G'_n\}^\infty_{n=1}$ even
the largest components have at most $K$ vertices.
We call a vertex $z\in V(G_n)=V(G_n')$ {\bf exceptional} if we removed at
least one of the edges incident to $z$. Clearly, the number of exceptional
vertices in $G_n$ is not greater than $\delta_1 d |V(G_n)|\,.$ Let $H_1,
H_2,\dots, H_{M_K}$ be the isomorphism classes of connected $(X,S)$-graphs with
vertex degree bound $d$ having at most $K$ vertices.
For $1\leq i \leq M_K$ we denote by $c^n_i$ the number of components
in $G_n'$ colored-isomorphic to $H_i$. Let $r^n_i=\frac{c^n_i}{|V(G_n)|}$ for
any $1\leq i \leq M_K$.
Pick a subsequence $\{G_{n_k}\}^\infty_{k=1}$ such that
for any  $1\leq i \leq M_K$ and $l,m\geq 1$:
\begin{equation} \label{e17}
|r_i^{n_l}-r_i^{n_m}|\leq \delta_2=\frac{\e}{100 M_K K}\,.
\end{equation}
Now pick two numbers $l,m\geq 1$. The proposition shall follow from the
following lemma.
\begin{lemma}
$\dr(G_{n_l},G_{n_m})\leq \e\,.$
\end{lemma}
\proof
First consider the graph $Z_l$ consisting of $|V(G_{n_m})|$ disjoint
copies of $G_{n_l}$ and the graph  $Z_m$ consisting of $|V(G_{n_l})|$ disjoint
copies of $G_{n_m}$. We assume that the vertex spaces of $Z_l$ and $Z_m$ are
the same. Also, we consider the subgraphs $Z_l'$ (resp. $Z_m'$) consisting
of $|V(G_{n_m})|$ (resp. $|V(G_{n_l})|$  copies of $G'_{n_l}$ (resp. 
$G'_{n_m})$. In $Z_l'$ (resp. in $Z_m'$) we have $|V(G_{n_m})| c_i^{n_l}$
(resp. $|V(G_{n_l})| c_i^{n_m}$) components isomorphic to $H_i$. 
For $1\leq i \leq M_K$ let $q_i=\min\{|V(G_{n_m})| c_i^{n_l}, 
|V(G_{n_l})| c_i^{n_m}\}$. Choose $q_i$ copies of components
of $Z'_l$ (resp. of $Z'_m)$ isomorphic to $H_i$. \\ If $z\in V(Z_l)$ 
 (resp.
$z\in V(Z_m))$ is {\it not} in the chosen copies for any $1\leq i 
\leq M_K$, then call $z$ (resp. $w$) a {\bf non-matching} vertex. Now construct
a permutation $\sigma$ on the vertices of $Z_m$ that for each 
$1\leq i \leq M_K$
maps the chosen $q_i$ copies of $Z'_m$ onto the $q_i$ copies of $Z'_l$ 
isomorphically. Define $\sigma$ arbitrarily on the non-matching vertices.

\noindent
{\bf Claim:}
\begin{equation}
\label{e18}
\delta(Z_l, Z^\sigma_m)\leq \frac{|A| + |B| + |C| +D|}{ |V(Z_l)|}\,,
\end{equation}
{\it where $A$ (resp. $B$) is the set of non-exceptional 
vertices in $Z_l$ (resp.
in $Z_m$) and $C$ (resp. $D$ is the set of non-matching vertices in 
 $Z_l$ (resp.
in $Z_m$).}

\noindent
Indeed, if $z$ is not in $A\cup \sigma(B)\cup C\cup \sigma (D)$ then
its star in $Z_l$ is the same as its star in $Z_m^\sigma$. 

\noindent
Now by our earlier observation
$$\max\{|A|,|B|\}\leq \delta_1 d |V(G_{n_l})| |V(G_{n_m})|=\delta_1 d 
|V(Z_l)|\,.$$
For the number of non-matching vertices we have the estimates
$$|C|\leq\sum^{M_K}_{i=1}\left(|V(G_{n_m})| c_i^{n_l}-q_i\right) K\,\quad
|D|\leq\sum^{M_K}_{i=1}\left(|V(G_{n_l})| c_i^{n_m}-q_i\right) K\,.$$
Recall that
$$\left|
  \frac{c^{n_l}_i}{|V(G_{n_l})|}-
\frac{c^{n_m}_i}{|V(G_{n_m})|}\right|\leq\delta_2\,.$$
That is $||V(G_{n_m})|c^{n_l}_i-|V(G_{n_l})|c^{n_m}_i|\leq 
\delta_2 |V(Z_l)|\,.$
Hence $|V(G_{n_m})| c_i^{n_l}-q_i\leq \delta_2 |V(Z_l)|\,.$
Therefore,
$$\frac{|C|}{|V(Z_l)|}\leq M_K\delta_2 K\,\quad\mbox{and}\quad
\frac{|D|}{|V(Z_m)|}\leq M_K\delta_2 K\,.$$
Thus by (\ref{e18})
$$\dr(G_{n_l}, G_{n_m})\leq \ds(Z_l,Z^\sigma_m)\leq 2d\delta_1+
2m\delta_2K\leq \e\,.$$
Hence our proposition follows. \qed

\noindent
We finish this subsection with a conjecture.
\begin{conjecture}
\label{conj1}
Let $\bG=\{G_n\}^\infty_{n=1}$ be a weakly convergent antiexpander sequence
of uniformly bounded vertex degrees. Then $\bG$ is
strongly convergent as well.
\end{conjecture}

\noindent
For the first sight 
the conjecture might seem to be a little bit too bold, nevertheless
it is not very hard to check that Proposition \ref{choose} together
with the Ornstein-Weiss Quasi-Tiling Lemma implies the conjecture for
the \Fo-subsets of Cayley-graphs of amenable groups.
\subsection{Examples of antiexpanders} \label{anti}

\noindent
\underline{\bf Graphs with subexponential growth:}\,
Recall that a monotone function $f:\bN\to \bN$ has {\bf subexponential growth}
if for any $\beta>0$ there exists $r_\beta>0$ such that $f(r)\leq exp(\beta
r)$
if $r\geq r_\beta\,.$
We say that a graph $G$ has growth bounded by $f$ if for any $x\in V(G)$ and
$r\geq 1$, $|B_r(x)|\leq f(r)$.
A graph sequence $\bG=\{G_n\}^\infty_{n=1}$ is of subexponential growth if
there exists a monotone function $f:\bN\to \bN$ of subexponential growth
such that $G_n$ has growth bounded by $f$ for any $n\geq 1$. The following
proposition is due to Jacob Fox and J\'anos Pach.
\begin{proposition}
If $\bG=\{G_n\}^\infty_{n=1}$ is a graph sequence of subexponential growth
with vertex degree bound $d$ then it is an antiexpander sequence.
\end{proposition} Fix a constant $\epsilon>0$.
\begin{lemma}
For any function $f$ of subexponential growth 
 there exists $R_f>0$ such that
if a graph $G$ has growth bounded by $f$ then for some $1\leq r \leq R_f$
$$\frac{|B_{r+1}(x)|}{|B_r(x)|}\leq 1+\epsilon\,.$$ \end{lemma}
\proof
If the lemma does not hold, then for each $R>0$ there exists a graph $G$ of 
growth bounded by $f$ such that for some $x\in V(G)$ :
$|B_R(x)|\geq(1+\epsilon)^{R-1}$
\, which is in contradiction with the subexponential growth condition. \qed

\noindent
By the vertex degree condition, it is enough to prove that if $G$ is an 
arbitrary finite graph of growth bounded by $f$ and of vertex degree
bound $d$ then there exists $T\subseteq V(G)$, $|T|\leq \epsilon |V(G)|$
such that the spanned subgraph of $V(G)\backslash T$ consists of components
containing at most $f(R_f)$ vertices. We use a simple induction.
If $|V(G)|=1$ the statement trivially holds. Suppose that the statement
holds for $1\leq i \leq n$ and let $|V(G)|=n+1$. Consider a vertex $x\in
V(G)$. Then there exists $1\leq r \leq R_f$ such that
$\frac{|B_{r+1}(x)\backslash B_r(x)|} {|B_r(x)|}\leq \epsilon$.
Consider $V(G)$ as the disjoint union
$$V(G)=B_r(x)\cup (B_{r+1}(x)\backslash B_r(x))\cup Z\,.$$
Obviously the vertices in $Z$ have no adjacent vertex in $B_r(x)$. Let $H$ be
the not necessarily connected subgraph spanned by the vertices of $Z$.
By induction, we have a set $S\subset Z$, $|S|\leq \epsilon |Z|$ such
that the components of the spanned subgraphs of $Z\backslash S$ 
are containing at most $f(R_f)$
vertices. Now let $T=(B_{r+1}(x)\backslash B_r(x))\cup S$. The components
of the graph spanned by $V(G)\backslash T$ are the components considered
above and $B_r(x)$. Clearly, $|B_r(x)|\leq f(R_f)$, hence the proposition
follows. \qed

\vskip0.1in
\noindent
\underline{\bf Amenability}
In \cite{Elekcost} we proved that if $G$
is a finitely generated residually finite group with a nested sequence
of finite index subgroups $\Gamma_n$; $\cap^\infty_{n=1}\Gamma_n$, then
their Cayley-graphs (with respect to generators of $\Gamma$) form an
antiexpander sequence if and only if $\Gamma$ is amenable. By the same
simple application of the Ornstein-Weiss Quasi-Tiling Lemma \cite{OW}
one can easily prove that if $\{G_n\}^\infty_{n=1}$ is a Folner-sequence
in the Cayley graph of a finitely generated amenable group $\Gamma$ then they
 form a weakly convergent
antiexpander sequence. Note that these graph sequences might have exponential
growth.

\section{Operators on graph sequences} \label{algebra}
\subsection{The weak convergence of operators}
Recall that for a finite graph $G$ the Laplacian $\Delta_G:L^2(V(G))\to
L^2((V(G))$ is a linear operator acting the following way
$$(\Delta_G f)(x)=\deg(x) f(x)-\sum_{x\sim y} f(y)\,.$$
In general, let us consider linear operators $A$ on the vertex set of a finite
a graph $G$ given by operator kernels
$A:V(G)\times V(G)\to \bR$:
$$Af(x)=\sum_{y\in V(G)} A(x,y)f(y)\,.$$
Note that we shall slightly abuse the notation and use the same letter for
an operator and its operator kernel.
Now let $\alpha\in\cA(d,r,X,S)$ and let $H$ be a rooted $(X,S)$-colored
graph representing the class $\alpha$. Suppose that $f:V(H)\to\bR$ is
a function that is invariant under the rooted colored automorphisms of 
$H$ (i.e.
$f(x)=f(\sigma(x))$ if $\sigma$ is a rooted colored automorphism).
Then the function $f$ can be viewed as a function on the vertices of $\alpha$.
We call such a function an {\bf invariant function} on $\alpha$.
Now let $\bG=\{G_n\}^\infty_{n=1}$ be a weakly convergent sequence of
$(X,S)$-graphs
of vertex degree bound $d$ and let $\bA=A_n:V(G_n)\times V(G_n)\to\bR$ be a
sequence
of operator kernels such that:
\begin{itemize}
\item
There exists a uniform bound $m_{\bA}>0$ such that $\sup_{n\in\bN} \sup_{x,y\in
  V(G_n)} |A_n(x,y)|\leq m_{\bA}\,.$
\item
There exists a uniform bound $s_{\bA}>0$ such that for any $n\geq 1$,
$x,y\in  V(G_n)$, $A_n(x,y)=0$ if $d_{G_n}(x,y)>s_{\bA}$.
\item For the same constant $s_{\bA}$, if $\alpha\in\cA(d,s_{\bA},X,S)$
and $\limn p_{G_n}(\alpha)\neq 0$, then there exists an invariant function
$f^{\bA}_\alpha$ on $\alpha$ 
such that if the $s_{\bA}$-neighborhood of a vertex
$x$ represents the class $\alpha$ then $A_n(x,.)$ represents 
$f^{\bA}_\alpha$ on
$B_{s_{\bA}}(x)$. That is the operators depend only on the local patterns.
\end{itemize}
Then we call the sequence $\bA=\{A_n\}^\infty_{n=1}$ be {\bf weakly convergent
operator sequence} on $\bG=\{G_n\}^\infty_{n=1}$.

\vskip0.1in
\noindent
\begin{example}
Let $\bG=\{G_n\}^\infty_{n=1}$ be a weakly convergent sequence of finite graphs
of bounded vertex degrees, then $\{\Delta_{G_n}\}^\infty_{n=1}$ is a weakly
convergent sequence of operators.
\end{example}
\begin{lemma}
Let $\bG=\{G_n\}^\infty_{n=1}$ be as above and $\bA=\{A_n\}^\infty_{n=1}$ and
$\bB=\{B_n\}^\infty_{n=1}$ be weakly convergent operator sequences.
Then both $\bA\bB=\{A_nB_n\}^\infty_{n=1}$ and 
$\bA+\bB=\{A_n+B_n\}^\infty_{n=1}$ 
are weakly convergent operator sequences. Also, $\bA^*=\{A^*_n\}^\infty_{n=1}$
is a weakly convergent operator sequence.
\end{lemma}
\proof
Note that
\begin{equation}\label{product}
A_nB_n(x,y)=\sum_{z\in V(G_n)} A_n(x,z)B_n(z,y)\,.
\end{equation}
Hence
\begin{itemize}
\item $A_nB_n(x,y)=0$ if $d_{G_n}(x,y)>s_{\bA}+s_{\bB}$.
\item $|A_nB_n(x,y)|\leq m_{\bA}m_{\bB} t$, where
$t$ is the maximal possible number of vertices of a ball of
 radius $s_{\bA}$ in a graph $G$ of vertex degree bound $d$.
\end{itemize}
Now let $y,z\in V(G_n)$ such that $B_{s_A+s_B}(y)\cong B_{s_A+s_B}(z)$
and these balls are represented by $\alpha$ so that $\limn
p_{G_n}(\alpha)\neq 0$. By the pattern invariance assumption and the equation
(\ref{product}),
$A_nB_n(y,.)$ and $A_nB_n(z,.)$ represent the same invariant function. This
shows that $\{A_nB_n\}^\infty_{n=1}$ is a weakly convergent sequence
of operators. The case of $\{A_n+B_n\}^\infty_{n=1}$ can be handled
the same way. Now let us turn to the adjoint sequence. Clearly,
$A_n^*(x,y)=A_n(y,x)$. Again, suppose that the $2s_{\bA}$-neighborhoods of 
$y,z\in
V(G_n)$ are rooted colored isomorphic and represent a class $\alpha$
such that $\limn
p_{G_n}(\alpha)\neq 0$.
Then it is easy to check that $A^*(y,.)$ and $A^*(z,.)$ represent the
same invariant function that is $\{A^*_n\}^\infty_{n=1}$ is a weakly
convergent
operator sequence. \qed

\vskip0.1in
By our previous lemma if $\bG=\{G_n\}^\infty_{n=1}$ is a weakly convergent
sequence of $(X,S)$ graphs with vertex degree bound $d$ then the weakly
convergent operator sequences form a unital $*$-algebra $\cP_{\bG}$.
\subsection{The trace}
Let $\bG=\{G_n\}^\infty_{n=1}$ be a weakly convergent sequence of
$(X,S)$-graphs
and $P_{\bG}$ be the $*$-algebra of weakly convergent operator sequences. We
call $\bB=\{B_n\}^\infty_{n=1}$ a nulloperator if for any
$\alpha\in\cA(d,s_{\bB},X,S)$ the associated invariant function 
$f^{\bB}_{\alpha}$ is
zero. Clearly, the nulloperators are exactly those weakly convergent operator
sequences, where
$$\limn\frac {|\{ (x,y)\mid B_n(x,y)\neq 0\}|}{|V(G_n)|}=0\,.$$
It is easy to see that the nulloperators form an ideal $\cN_{\bG}$ in
$\cP_{\bG}$.
For a weakly convergent operator sequence $\bA=\{A_n\}^\infty_{n=1}$ we
define 
$$\trg (\bA):=\limn\frac{1} {|V(G_n)|}\sum_{x\in V(G_n)} A_n(x,x)=\limn
\frac{1} {|V(G_n)|} \tr(A_n)\,.$$
\begin{proposition}
\begin{enumerate}
\item The limit in the definition of $\trg(\bA)$ does exist.
\item $\trg$ is a {\bf trace}, that is a linear functional satisfying
$\trg(\bA\bB)=\trg(\bB\bA)$\,.
\item $\trg$ is faithful on the quotient space $\cP_{\bG}/\cN_{\bG}$
that is $\trg(\bA^*\bA)>0$ if $\bA\notin \cN_{\bG}$ and
$\trg(\bB)=0$ if $\bB\in \cN_{\bG}$.
\end{enumerate}
\end{proposition}
\proof
Let $\alpha\in\cA(d,s_{\bA},X,S)$ and $T(G_n,\alpha)$ be the set of vertices
in $V(G_n)$ such that $B_{s_{\bA}}(x)$ is represented by $\alpha$.
$$\frac {1} {|V(G_n)|}\sum_{x\in V(G_n)} A_n(x,x)=
\frac {1} {|V(G_n)|}\sum_{\alpha\in\cA(d,s_{\bA},X,S)} (\sum_{x\in T(G_n,\bA)}
A_n(x,x))=
$$ $$=
\sum_{\{\alpha\in\cA(d,s_{\bA},X,S)\,\mid\, \limn p_{G_n}(\alpha)\neq 0\}}
p_{G_n}(\alpha)f^{\bA}_\alpha(r) +o(1)\,,$$
where $f^{\bA}_\alpha(r) $ is the value of the invariant function at the root.
Thus the limit exists and equals to
\begin{equation}
\label{trace}
\sum_{\{\alpha\in\cA(d,s_{\bA},X,S)\,\mid\, \limn p_{G_n}(\alpha)\neq 0\}}
p_{G_n}(\alpha)f^{\bA}_\alpha(r)\,. \end{equation}

Observe that second statement of the proposition
 follows immediately from the fact that 
$\tr(A_n B_n)=\tr(B_n A_n)$.
\vskip0.1in
\noindent
Let $\bB\in\cN_{\bG}$ then clearly
$$\limn \frac {1} {|V(G_n)|}\sum_{x\in V(G_n)} B_n(x,x)=0 $$
hence $\trg(\bB)=0$.
Now let $\bB\notin\cN_{\bG}$.
Then
$$\frac {1} {|V(G_n)|}\tr(B_n^*B_n)= \frac {1} {|V(G_n)|}\tr(B_nB_n^*)=
\frac {1} {|V(G_n)|} 
\sum_{x\in V(G_n)} \sum_{y\in V(G_n)} B_n(x,y) B_n^*(x,y)=$$
$$=\frac {1} {|V(G_n)|} \sum_{x\in V(G_n)} (\sum_{y\in V(G_n)} B_n(x,y)^2)\,.$$
If there exists at least one $\alpha\in \cA(d,s_{\bB},X,S)$ such that
$\limn p_{G_n}(\alpha)\neq 0$ and the associated invariant function
$f^{\bB}_\alpha$ is non-zero then
$$\limn \frac {1} {|V(G_n)|}\tr(B_n^*B_n)\geq \limn p_{G_n}(\alpha)
\sum_{v\in V(\alpha)} |f^{\bB}_\alpha (v)|^2>0\,.$$

\subsection{The von Neumann algebra of a graph sequence}
The von Neumann algebra of $\bG$, $N_{\bG}$ is constructed by the
GNS-construction the usual way. The algebra $\cP_{\bG}/\cN_{\bG}$ is
a pre-Hilbert space with inner product
$$\langle [\bA], [\bB]\rangle =\trg(\bB^* \bA)\,,$$
where $[\bA]$ denotes the class of $\bA$ in $\cP_{\bG}/\cN_{\bG}$\,.
Then $L_{[\bA]}[\bB]=[\bA\bB]$ defines a representation of
$\cP_{\bG}/\cN_{\bG}$ on this pre-Hilbert space.
\begin{lemma}
$L_{[\bA]}$ is a bounded operator for any $\bA\in\cP_{\bG}$.
\end{lemma}
\proof
Let $\bA,\bB\in P_G$. Denote by $\|\,\|$ the pre-Hilbert space norm. Then
\begin{equation} \label{l32}
\|\bB\|^2=\tr_G(\bB^*\bB)=\trg(\bB\bB^*)=\limn \frac{1}{|V(G_n)|}
\sum_{x\in V(G_n)}(\sum_{y\in V}
|\bB(x,y)|^2)\,.\end{equation}
$$\|L_\bA \bB\|^2=\langle \bA\bB, \bA\bB \rangle =
\trg((\bA\bB)^*\bA\bB)=\trg(\bB^*\bA^*\bA\bB)= 
\trg(\bB\bB^*\bA^*\bA)=$$ $$=
\limn \frac{1}{|V(G_n)|}\sum_{x\in V(G_n)}(\sum_{y\in V(G_n)}
B_nB_n^*(x,y)A_n^*A_n(y,x))\leq $$ $$\leq\limn \frac{m_{\bA^*\bA}}{|V(G_n)|}
\sum_{x\in V(G_n)}
(\sum_{\{y\in V(G_n)\,\mid\, d_{G_n}(x,y)\leq s_{\bB\bB^*}\}}
|B_nB_n^*(x,y)|)\,.$$
Note that
$$|B_nB_n^*(x,y)|=|\sum_{z\in V(G_n)} B_n(x,z) B_n^*(z,y)|\leq
\sum_{z\in V(G_n)}|B_n(x,z)| |B_n(y,z)|\leq $$ $$ \leq\frac{1}{2}
\sum_{z\in V(G_n)}(|B_n(x,z)|^2+|B_n(y,z)|^2)\,.$$
Let $t^n_x:=\sum_{y\in V(G_n)} |B_n(x,y)|^2\,.$ Then by (\ref{l32})
$$\|\bB\|^2=\limn \frac{1}{|V(G_n)|}\sum_{x\in V(G_n)}t^n_x\,.$$
Hence
$$\|L_\bA \bB\|^2\leq \limn \frac{m_{\bA^*\bA}}{|V(G_n)|}
\sum_{x\in V(G_n)}(\sum_{\{y\in V(G_n)\,\mid\,
 d_{G_n}(x,y)\leq s_{\bB\bB^*}\}}\frac{1}{2}
(t^n_x+t^n_y))\leq $$$$ \leq \limn \frac{m_{\bA^*\bA}T} {|V(G_n)|}
\sum_{x\in V(G_n)}t^n_x\,,$$
where $T$ is the maximal possible number of vertices in an 
ball of radius $s_{\bB\bB^*}$ in a graph of vertex degree bound $d$.
 Hence $\|L_{[\bA]} [\bB]\|^2\leq m_{\bA^*\bA}T \|\bB\|^2\,.$\qed

\vskip0.1in
\noindent
Thus
$\cP_{\bG}/\cN_{\bG}$ is represented by bounded operators on the Hilbert-space
closure. Now the von Neumann algebra $N_{\bG}$ is defined as the weak-closure
of $\cP_{\bG}/\cN_{\bG}$. Then the trace 
$$\trg ([\bA])=\langle L_{[\bA]} {\bf 1}, {\bf 1}\rangle$$ extends to
$N_{\bG}$ as an ultraweakly continuous, faithful trace. 
If $\bA=\{A_n\}^\infty_{n=1}$ is a weakly convergent operator sequence
then we call $[\bA]\in N_{\bG}$ the {\bf limit operator} of $\bA$.
\subsection{Representation in the von Neumann algebra and the integrated 
density of states}
In the whole subsection let $\bG=\{G_n\}^\infty_{n=1}$ be a weakly convergent
subsequence of finite connected graphs with vertex degree bound $d$ and 
$\bB=\{B_n\}^\infty_{n=1}$ be a weakly convergent sequence of self-adjoint
operators.
\begin{lemma}
$[\bB]\in N_{\bG}$ is a self-adjoint
operator.
\end{lemma}
\proof
$$\langle L_{[B]}[X],[Y]\rangle=\limn \frac {1} {|V(G_n)|} \tr(Y_n^*B_nX_n)
$$
$$\langle [X], L_{[B]}[Y]\rangle=\limn \frac {1} {|V(G_n)|}
\tr((B_nY_n)^*X_n)=\limn \frac {1} {|V(G_n)|} \tr(Y_n^*B_nX_n)\quad\qed$$
\begin{lemma}
There exists some $K>0$ depending on $\bB$ such that
$\|B_n\|\leq K$ for any $n\geq 1$.
\end{lemma}
For the unit vectors $f,g\in L^2(V(G_n))$
$$|\langle B_n(f),g\rangle| =|\sum_{x,y\in V(G_n)} B_n(x,y)f(x)g(y)|\leq
m_{\bB} |\sum_{\{x,y\in V(G_n)\,\mid\, d_{G_n}(x,y)\leq 
s_{\bB}\}} f(x)g(y)|\leq$$
$$\leq m_{\bB}\, |\sum_{\{x,y\in V(G_n), d_{G_n}(x,y)
\leq s_{\bB}\}} f^2(x)+g^2(y)|\,.$$

The number of occurrences of $f^2(x)$ resp. $g^2(y)$ on the right hand side
is not greater than $t(d,s_{\bB})$ the maximal possible number of vertices
in a ball of radius $s_{\bB}$ of in a graph of vertex degree bound $d$.
Hence
$$\|B_n\|\leq 2 m_{\bB}t(d,s_{\bB})\,.\quad\qed$$

Since $[\bB]\in N_{\bG}$ is a self-adjoint element of a von Neumann
algebra we have the spectral decomposition:
$$[\bB]=\int^\infty_0\lambda\,dE^{[\bB]}_\lambda\,,$$
where $E^{[\bB]}_\lambda\in N_{\bG}$ is the associated spectral projections
$\chi_{[0,\lambda]}([\bB])$. The {\bf spectral measure} $\mu_{[\bB]}$
is defined by
$$\mu_{[\bB]}([0,\lambda]):=\trg(E^{[\bB]}_\lambda)\,.$$
Note that the spectral distribution of $\mu_{[\bB]}$ is just
$N_{[\bB]}(\lambda)=\mu_{[\bB]}([0,\lambda])\,.$ Let us notice
that the spectral distributions $N_{B_n}$ also define probability measures
on $\bR$ by $\mu_{B_n}[-\infty,\lambda]= N_{B_n}(\lambda)\,.$
\begin{theorem}
\label{integrated}
(the existence of the integrated density of states)
The probability measures $\mu_{B_n}$ weakly converge to $\mu_{[\bB]}$.
Thus for any continuity point $\lambda$ of $N_{[\bB]}$ we have
the Pastur-Shubin formulas
$$\limn N_{B_n}(\lambda)=\trg (E^{[\bB]}_\lambda)\,.$$\end{theorem}
\proof
If $P\in\bR[x]$ is a polynomial then
$$\int^a_{-a} P(\lambda)\, d\mu_{B_n}(\lambda)=\frac {1} 
{|V(G_n)|}\tr(P(B_n))\,$$
and
$$\int^a_{-a} P(\lambda)\, d\mu_{\bB}(\lambda)=\frac {1} 
{|V(G_n)|}\trg(P(\bB))\,,$$
where 
$$a:=\max\{\|[\bB]\|,\sup_{n\geq 1} \|B_n\|\}\,.$$
By definition, $\limn \frac {1} {|V(G_n)|} \tr(B^k_n)=\trg([B]^k)$ thus
$$\limn \int^a_{-a} P(\lambda)\, d\mu_{B_n}(\lambda)= 
\int^a_{-a} P(\lambda)\, d\mu_{\bB}(\lambda)\,.$$
That is $\{\mu_{B_n}\}^\infty_{n=1}$ weakly converge to $\mu_{\bB}$. \qed
\subsection{Strong convergence of graphs implies the uniform convergence
of the spectral distribution functions}
The goal of this subsection is to prove the following proposition.
\begin{proposition}
Let $\bG=\{G_n\}^\infty_{n=1}$ be a strongly convergent graph sequence of
vertex degree bound $d$. Let $\bB=\{B_n\}^\infty_{n=1}$ be a weakly
convergent sequence of operators on $\bG=\{G_n\}^\infty_{n=1}$. Then the
spectral distribution functions $N_{B_n}$ uniformly converge.
\end{proposition}
\proof First we need an elementary lemma on small rank perturbations.
\begin{lemma}\label{rank}
Let $C:\bR^n\to\bR^n$, $D:\bR^n\to\bR^n$ be self-adjoint linear
transformations
such that $\rk(C-D)\leq \e n$. Let $N_C$ resp. $N_D$ be the normalized
spectral distribution functions of $C$ resp. of $D$. Then
$$\|N_C-N_D\|_{\infty}\leq \e\,.$$ \end{lemma}
\proof
Since $N_C(\lambda)=N_{C+r\,Id}(\lambda+r)$ we can suppose that both $C$
and $D$ are positive operators. We call $V\subseteq \bR^n$ an $M_C$-subspace
if for any $v\in V$: $\|Cv\|\leq \lambda \|v\|\,.$
Let $d_C(\lambda)$ denote the maximal dimension of a $M_C$-subspace.
Observe that $d_C(\lambda)=n N_C(\lambda)$ that is the number of eigenvalues
not greater than $\lambda$. Indeed, if $W^+_\lambda$ is the space spanned
by the eigenvectors belonging to eigenvalues larger than $\lambda$ then
$V\cap W^+_\lambda=0$ for any $M_C$-subspace $V$. Therefore 
$\dim_V\leq n-\dim W^+_\lambda\,.$
On the other hand, $\dim W^-_\lambda=n- \dim W^+_\lambda\,,$ where
$W^-_\lambda$ is the subspace spanned by the eigenvectors belonging to
eigenvalues not greater than $\lambda$. Since $W^-_\lambda$ is
an $M_C$-subspace $d_C(\lambda)=n N_C(\lambda)$.

\noindent
If $\rk(C-D)\leq \e n$, then $\dim \ker (C-D)\geq n-\e n$.
Let $Q_\lambda=W^-_\lambda\cap \ker (C-D)$. Then 
$\dim Q_\lambda\geq \dim W^-_\lambda-\e n$ and $Q_\lambda$ is
an $M_D$-subspace of $\bR^n$. Therefore
$|N_C(\lambda)-N_D(\lambda)|\leq\e$.\quad\qed

\noindent
Now let us return to the proof of our proposition.
Suppose that $\dr(G_n,G_m)<\frac{\e}{C}$, where $C$ is the maximal
possible number of vertices of a ball of radius $s_{\bA}$ in a graph of
vertex degree bound $d$.
Let $q,r$ be natural numbers such that $\ds(qG_n,rG_m)<\frac{\e}{C}$.
We define the linear operator $\ah_n$ the following way.
\begin{itemize}
\item $\ah_n(\hat{x},\hat{y})=A_n(x,y)$ if
$\hat{x}$ and $\hat{y}$ are the vertices in a component 
 of $qG_n$ corresponding to
the vertices $x$ and $y$.
\item $\ah_n(\hat{z},\hat{w})=0$ if $\hat{z}$ and $\hat{w}$ are not in the
  same component.
\end{itemize}
Clearly, the normalized spectral distribution of $\ah_n$ is exactly the
same as the one of $A_n$. We define $\ah_m$ similarly on $rG_m$. Now
let $\ah_m^\sigma(\hat{z},\hat{w})$ defined as
$\ah_m(\sigma(\hat{z}),\sigma(\hat{w}))$.
\begin{lemma}
\label{kulcs}
Suppose that $\delta(qG_n,(rG_m)^\sigma)<\frac{\e}{C}$, then
$\rk(\ha_n-\ha_m^\sigma)<\e |V(qG_n)|$.
\end{lemma}
\proof
It is enough to prove that
$$\dim\ker (\ha_n-\ha_m^\sigma)>|V(qG_n)|-\e |V(qG_n)|\,.$$
Let $x\in V(qG_n)$ and $v_x$ be the vector in $L^2(qG_n)$ having
the value $1$ at $x$ and $0$ everywhere else. Clearly, if the
$s_{\bA}$-neighborhood of $x$ in $qG_n$ is the same as the $s_{\bA}$
neighborhood of $x$ in $(rG_m)^\sigma$ then
$(\ha_n-\ha_m^\sigma)v_x=0$. If $x$ is not such a vertex, then $x$ in the
$s_{\bA}$-neighborhood of an exceptional vertex of $qG_n$. Thus the number
of such vertices is not greater than $C \frac{\e}{C} |V(qG_n)|\,.$
Thus the lemma follows. \qed

\noindent
By Lemma \ref{rank}, the normalized spectral distributions
$\{N_{A_n}\}^\infty_{n=1}$ converge,
that proves our proposition. \qed

\section{Finite range operators on infinite graphs} \label{infinite}
\subsection{Antiexpander graphs}

We say that an amenable (see Introduction)
 $(X,S)$-graph has {\bf uniform patch frequency} (UPF)
if all of its \Fo subgraph sequences are weakly convergent.
We call $G$ an an UPF-antiexpander if :
\begin{itemize}
\item
$G$ is an amenable graph with the $UPF$-property.
\item
All \Fo\, subgraph sequences are antiexpanders.
\end{itemize}
\begin{example}
Let $\Gamma$ be a finitely generated amenable group and \\
$S=\{s_1,s_2,\dots,s_{d},s_1^{-1}, s_2^{-1},\dots, s_d^{-1}\}$ be a symmetric
generating system. Consider the Cayley-graph Cay$(\Gamma,S)$, where
\begin{itemize}
\item $V(\mbox{Cay}(\Gamma,S))=\Gamma\,.$
\item $(x,y)\in E(\mbox{Cay}(\Gamma,S))$ if $\gamma x=y$ for some
$\gamma\in S$. \\ In this case the edge color of $(x,y)$ is $\gamma$.
\end{itemize}
Then $\mbox{Cay}(\Gamma,S))$ is an UPF-antiexpander (see the the end of
Subsection \ref{anti}).
\end{example}
\begin{example}
The Delone-graphs in \cite{Lenz},\cite{LS} are antiexpanders (since they
have polynomial \\ growth) and they have the UPF-property.
\end{example}
\subsection{Self-similar graphs}
The goal of this subsection is to provide ample amount of examples of
UPF-antiexpanders, namely {\bf self-similar graphs} (see e.g. \cite{Kron} and
the references therein for similar constructions).
 First we fix two positive integers $d$ and $k$.
Let $G_1$ be a finite connected graph with
vertex degree bound $d$ with a distinguished subset
$S_1\subset V_1(G_1)$, which we call the set of connecting vertices.
Now we consider the graph $\widetilde{G}_1$, which consists of $k$
disjoint copies of $G_1$ with following additional properties:
\begin{itemize}
\item The graph $G_1$ is identified with the first copy.

\item In each copy the vertices associated to a connecting vertex
of $G_1$ is a connecting vertex of the graph $\widetilde{G}_1$.
\end{itemize}
The graph $G_2$ is defined by adding some edges to $\widetilde{G}_1$
such that both endpoints of these new edges are connecting vertices.
We must also ensure that the resulting graph  still has
 vertex degree bound $d$.
Finally the subset $S_2\subset V(G_2)$ is chosen as a subset of the
connecting vertices of $\widetilde{G}_1$ such that $S_2\cap
V(G_1)=\emptyset.$ That is $G_1\subset G_2$ is a subgraph and the connecting
vertices of $G_2$  are not in $G_1$.
Inductively, suppose that the finite graphs $G_1\subset G_2\subset
\dots \subset G_n$ are already defined and the vertex degrees in
$G_n$ are not greater than $d$. Also suppose that a set $S_n\subset
V(G_n)$ is given and $S_n\cap V(G_{n-1})=\emptyset$. Now the graph
$\widetilde{G}_n$ consists of $k$ disjoint copies of $G_n$ and
\begin{itemize}
\item The graph $G_n$ is identified with the first copy.

\item In each copy the vertices associated to a connecting vertex
of $G_n$ is a connecting vertex of the graph $\widetilde{G}_n$.
\end{itemize}

Again, $G_{n+1}$ is constructed by adding edges to $\widetilde{G}_n$
with endpoints which are connecting vertices, preserving the vertex
degree bound condition. Then the set of connecting vertices $S_{n+1}$ is
chosen as a subset of the connecting vertices of $\widetilde{G}_n$,
such that $S_{n+1}\cap V(G_n)=\emptyset\,.$ The union of the graphs
$\{G_n\}^\infty_{n=1}$ is connected infinite graph with vertex
degrees not greater than $d$. We call the graph $G$ {\it self-similar} if
$\limn \frac{|S_n|}{|V(G_n)|}=0$. 
\begin{proposition}
Self-similar graphs are UPF-antiexpanders.
\end{proposition}
\proof
First of all note that $\partial G_n\subseteq S_n$, hence self-similar graphs
are amenable.
\begin{lemma}
For any $\e>0$ there exists $\delta_\e>0$ and $q_\e\in\bN$ such that if
$H\subset G$ is a finite spanned subgraph and $\frac{|\partial
  H|}{|V(H)|}<\delta_\e$
then $\dr(H, G_{q_\e})<\e$. \end{lemma}
\proof
By our construction, 
for any $n\geq 1$, $V(G)=\cup^\infty_{i=1} V(G^i_n)$, where
\begin{itemize}
\item $G_n^i$ (as a spanned subgraph) is isomorphic to $G_n$.
\item $\frac{|\partial
  G^i_n|}{|V(G^i_n)|}\leq c_n=\frac{|S_n|}{|V(G_n)|}\,.$

\end{itemize}
Now we choose $q_\e$ so large that $c_{q_\e}<\frac{\e}{2}\,.$
Then let $\delta_\e=\frac{\e}{10 |V(G_{q_\e})|}$. We suppose
that $\frac{|\partial H|}{|V(H)|}<\delta_\e$.

\noindent
Consider the set of indices $I_H$ such that $i\in I_H$ if 
and only if $V(G^i_{q_\e})\cap V(H)\neq\emptyset$.
Also, let $i\in J_H\subseteq I_H$ if and only if
$V(G^i_{q_\e})\subset V(H)$. If $i\in J_H$, then we call $G^i_{q_\e}$ an
{\bf inner copy} of $G_{q_\e}$ in $H$. Clearly,
\begin{equation} \label{becs}
|J_H| |V(G_{q_\e})|\leq |V(H)|\leq |I_H| |V(G^i_{q_\e})|\,.
\end{equation}

Note that if $i\in I_H\backslash J_H$ then there exists an element
of $x\in\partial H\cap V(G_{q_\e})$. Therefore $|I_H\backslash J_H|\leq
|\partial H|\,.$
By our assumption,
$$|I_H\backslash J_H| | V(G_{q_\e})|\leq \delta_\e | V(G_{q_\e})| |V(H)|\,.$$
Hence by (\ref{becs})
\begin{equation} \label{inner}
|V(H)|-|J_H|| V(G_{q_\e})|\leq \delta_\e | V(G_{q_\e})| |V(H)|\,.
\end{equation}
Now we estimate the geometric graph distance of $H$ and $G_{q_\e}$. First
consider $| V(G_{q_\e})|$ disjoint copies of $H$ and
$|V(H)|$ disjoint copies of $G_{q_\e}$ on the same vertex set $V$ of
cardinality $| V(G_{q_\e})| |V(H)|$.
Choose a permutation $\sigma\in S(V)$ that maps each of the 
$|J_H|| V(G_{q_\e})|$ inner copies of $G_{q_\e}$ in $| V(G_{q_\e})| H$
isomorphically into one of the $|V(H)|$ copies of $G_{q_\e}$.
Observe that if the star of a vertex $x\in V$ is not the same in
$| V(G_{q_\e})| H$ as in $|V(H)|G_{q_\e}$ then 
\begin{itemize}
\item
either $x$ is not in one of the inner copies
\item 
or $x$ is on the boundary of one of the inner copies.
\end{itemize}
Hence by (\ref{inner}) and our assumption
$$\delta\left( (|V(G_{q_\e})| H)^\sigma\right), |V(H)| G_{q_\e})\leq
\frac{\delta_\e |V(G_{q_\e})|^2 |V(H)|+ c_{q_\e} |V(G_{q_\e})||V(H)|}
{|V(G_{q_\e})||V(H)|}\,.$$
That is 
$$\dr(H, G_{q_\e})\leq \delta_\e |V(G_{q_\e})| + c_{q_\e}\leq \e\,.\quad\qed$$
Now we turn back to the proof of our proposition.
It is enough to prove that if $\{H_n\}^\infty_{n=1}$ is a \Fo-subgraph
sequence then it is Cauchy in the $\dr$-metric.
However by the previous lemma, if $n,m$ are large enough then
$\dr(H_n, G_{q_\e})\leq e$ and $\dr(H_m, G_{q_\e})\leq \e$, that is
$\dr(H_n,H_m)\leq 2\e$. \qed

\subsection{The main result}
Let $G$ be an infinite connected $(X,S)$-graph
with bounded vertex degrees and $A:V(G)\times V(G)\to \bR$ be
an operator kernel. We call $A$ a pattern-invariant finite range operator
if there exists some $s_A$ such that 
\begin{itemize}\item
$A(x,y)=0$ if $d_G(x,y)>s_A$
\item $A(x,y)=A(\phi(x),\phi(y))$ if $\phi$ is a rooted color isomorphism
from $B_{s_A}(x)$ to $B_{s_A}(\phi(x))$.
\end{itemize}
Let $\bG=\{G_n\}^\infty_{n=1}$ be spanned subgraph sequence in $G$. Then
the {\bf finite volume approximation} of $A$ on $\bG$ is
the sequence $A_n:V(G_n)\times V(G_n)\to \bR$, where
$A_n(x,y)=A(x,y)$. In \cite{LS} Lenz and Stollmann proved that if
$A$ is a pattern-invariant self-adjoint operator on a Delone-graph then
the spectral distribution functions of $\{A_n\}^\infty_{n=1}$ associated
to an arbitrary \Fo-sequence converge uniformly to an integrated density
of state.

\noindent
Let $s\in \bR\Gamma$ be a self-adjoint element of the group algebra of
a finitely generated amenable group $\Gamma$, then $s$ defines a
pattern-invariant self-adjoint operator $A^s$ on  Cay$(\Gamma,S)$ by
\begin{itemize}
\item $A^s(x,y)=c_{\gamma}$ if $y=\gamma x$ and $s=\sum c_\gamma \gamma$
(note that the self-adjointness of $s$ means that $c_\gamma=c_{\gamma^{-1}}$).
\end{itemize}
Again we have an associated sequence of finite dimensional approximation
for any Folner sequence $\{G_n\}^\infty_{n=1}$. 
Then the spectral distributions
of the operators $\{A^s_n\}^\infty_{n=1}$ are converge uniformly to the
spectral measure of $s$ in the von Neumann algebra of the group $\Gamma$
\cite{DLM}.
Our main result for infinite graphs generalizes the two results above
answering Question 1. for a large family of graphs.
\begin{theorem}
\label{main}
Let $G$ be an infinite antiexpander $(X,S)$ graph with uniform
patch-frequency. Let $A$ be a self-adjoint finite range pattern-invariant
operator on $G$. Then for any \Fo-sequence $\{G_n\}^\infty_{n=1}$ the
normalized spectral distributions of the finite volume approximations
$\bA=\{A_n\}^\infty_{n=1}$ uniformly converge to an integrated density of
state that does not depend on the choice of the \Fo-sequence.
\end{theorem}
\proof
Let $\bG=\{G_n\}^\infty_{n=1}$ and $\bG'=\{G'_n\}^\infty_{n=1}$ be two 
\Fo-sequences
in $G$. with associated approximating operator sequences
$\bA=\{A_n\}^\infty_{n=1}$ and $\bA'=\{A'_n\}^\infty_{n=1}$. Let
$X'=X\cup \{t\}$ and
let $H=\{H_n\}^\infty_{n=1}$ resp. $H'=\{H'_n\}^\infty_{n=1}$ be the
$(X',S)$-graph sequences obtained by recoloring the vertices in $\partial G_n$
resp. in  $\partial G'_n$ by the extra color $t$. Observe that $\bA$ resp.
$\bA'$ are weakly convergent operator sequences on $H$ resp. on $H'$.
Note that if $\alpha\in\cA(d,s_A,X',S)$ contains a vertex colored by $t$
then $\limn p_{H_n}(\alpha)=\limn p_{H'_n}(\alpha)=0$. Also, if $\limn
p_{H_n}(\alpha)\neq 0$ then $f^{\bA}_\alpha=f^A_\alpha$, where
$f^A_\alpha$ is the invariant function on $\alpha$ associated to the
finite range operator $A$.
\begin{lemma}
\label{coincide}
The spectral measures of the limit operator $[\bA]$
 of $\bA=\{A_n\}^\infty_{n=1}$ in
$\bN_H$ resp. of the limit operator $[\bA']$  of 
$\bA'=\{A'_n\}^\infty_{n=1}$ in
$\bN_{H'}$ coincide.
\end{lemma}
By Lemma \ref{trace},
$$\tr_{\bH}([\bA]^k)=\limn \frac{1} {|V(H_n)|} \tr(A_n^k)=
\sum_{\alpha\in\cA(d,ks_A,X',S)} \limn p_{H_n}(\alpha)
f_\alpha^{\bA^k}(r) \,,$$
where $f^{\bA^k}_\alpha (r)$ is the value at the root of the invariant
function on $\alpha$ associated to $\bA^k$.
Similarly,
$$\tr_{\bH'}([\bA']^k)=\limn \frac{1} {|V(H'_n)|} \tr((A'_n)^k)=
\sum_{\alpha\in\cA(d,ks_A,X',S)} \limn p_{H_n}(\alpha)f_\alpha^{\bA'^k}(r)\,.$$
By the pattern invariance of the finite range operator $A$ for any
$\alpha\in A(d,ks_A,X',S)$,
$$\limn p_{H_n}(\alpha)=\limn p_{H'_n}\quad\mbox{and}\quad
f_\alpha^{\bA^k}(r)=f_\alpha^{(\bA')^k}(r)\,.$$
Therefore the lemma follows from the definition of the spectral
measure. \qed

\vskip0.1in
Consequently, by Theorem \ref{integrated} we have the following corollary.
\begin{corollary}
Let $G$,$A$, $\bA=\{A_n\}^\infty_{n=1}$ be as in Theorem \ref{main}. Then
the normalized spectral distribution functions $\{N_{A_n}\}^\infty_{n=1}$
converge in any continuity point of a monotone right-continuous function
$g$ that does not depend on the choice of the \Fo-sequence.
\end{corollary}
Now suppose that $\{N_{A_n}\}^\infty_{n=1}$ does not converge uniformly
to $g$. Then there exists a subsequence $\{N_{A_{n_k}}\}^\infty_{k=1}$
such that for each $k\geq 1$,
$$\|N_{A_{n_k}}-g\|_\infty>\e>0\,.$$

 By our previous lemma and Proposition
\ref{choose} we can pick a subsequence of $\{N_{A_{n_k}}\}^\infty_{k=1}$
which uniformly converges. The limit function $f$ is right-continuous and
at the continuity points of $g$, $f$ must coincide with $g$. Thus $g=f$,
leading to a contradiction. This proves Theorem \ref{main}.\qed

\vskip0.1in
\noindent
{\bf Remark:} If we apply Theorem \ref{main} for Laplacian operators as in
Question 1 we
encounter a small difficulty. Namely, the finite volume approximation
operators of the Laplacian of $G$ 
on the \Fo-subgraphs $\{G_n\}^\infty_{n=1}$ are {\it not} 
the Laplacian operators 
$\Delta_{G_n}$
, since the degree
of a vertex in the subgraph and in the original graph are different if the
vertex is on the boundary.
Nevertheless we have the estimate 
$$\rk(p_n\Delta_G i_n-\Delta_{G_n})\leq |\partial{G_n}|\,.$$ That is by
our Lemma \ref{rank}, the limits of the normalized spectral distributions
of $\{p_n\Delta_G i_n\}^\infty_{n=1} $ and of $\{ \Delta_{G_n}\}^\infty_{n=1}$
coincide.

\end{document}